\documentclass[a4paper,12pt]{article}
\usepackage[english]{babel}

\title{Weak perturbations of shock waves.}
\author{Rasskazov I.O. \\Institute of Mathematics and Computer
Center,\\ Ufa, Russia, e-mail: rasskazov@imat.rb.ru}
\date{}
\begin{document}
\maketitle

\section*{Introduction.}

The Cauchy problem is considered for the perturbed strictly
hyperbolic 2x2 system of quasilinear equations
  $$
\left\{\begin{array}{ll}
u_t+\lambda(u)u_x=\varepsilon f(u,v), \
u|_{t=0}=\stackrel{o}{u}(x),
\\ v_t+\mu(u,v)v_x=\varepsilon g(u,v), \ v|_{t=0}=\stackrel{o}{v}(x),
 \ x\in R. \end{array}\right. \ \ \ 0<\varepsilon \ll 1 \eqno{(1)}$$
Here $\varepsilon$ is a small parameter. We assume, that the
initial data ${\stackrel{o}{u}}(x), \ {\stackrel{o}{v}}(x)$ have
jumps at $x=0$. The unperturbed problem (with $\varepsilon=0$) has
a persistent solution with two discontinuity lines (shock waves
\footnote{We will name the lines of discontinuity of the solution
in the plain $(x,t)$ as shock lines.}). Well-known Hugoniot
conditions are necessary for a uniqueness of the solution of the
problem (1). Both an asymptotics of shock waves position in the
plane (x,t) and an asymptotics of the perturbed problem solution
are discussed, when $\varepsilon\rightarrow 0$.
  \\
In this case the asymptotic formula
$$s^\pm(t,\varepsilon) \sim \sum_{i=0}^\infty \varepsilon^i
s^\pm_i(t), \ \varepsilon\rightarrow 0. \eqno{(2)}$$ holds for the
discontinuity lines, and
$$ u(x,t,\varepsilon) \sim \sum_{i=0}^\infty \varepsilon^i u_i(x,t),
 \ \ v(x,t,\varepsilon) \sim \sum_{i=0}^\infty \varepsilon^i v_i(x,t),
 \ \varepsilon\rightarrow 0\eqno{(3)}$$
hold for the solution in the continuity domain. Here the leading
terms $s^\pm_o(t)$,  $u_o(x,t)$ and $v_o(x,t)$ are the shock waves
and the solution of the unperturbed problem respectively.
\\
The corrections in the asymptotic solution $s^\pm_i(t), \
u_i(x,t), \ v_i(x,t)$ are found from linear equations.

We will use the following domains in the plain $(x,t)$
$$D^+=\{(x,t) \ | \ t>0, x>s^+(t,\varepsilon)\}, \quad
D^-=\{(x,t) \ | \ t>0, x<s^-(t,\varepsilon)\},$$
$$\widetilde D=\{(x,t) \ | \ t>0,
s^-(t,\varepsilon)<x<s^+(t,\varepsilon)\}.$$

\section*{Asymptotic solution in domains $D^{\pm}$ .}

To obtain an asymptotics in the domains $D^\pm$ we consider two
Cauchy problems for the system (1). The initial data are taken in
the form

(A) $u|_{t=0}=\stackrel{o}u(x), v|_{t=0}=\stackrel{o}v(x), \ x<0$

(B) $u|_{t=0}=\stackrel{o}u(x), v|_{t=0}=\stackrel{o}v(x), \ x>0$

Each of these problems can be solved (see [1]) in some area. The
Cauchy problem (1),(A) has a continuous solution in a domain
$\widetilde D^- $ in the plain $(x,t)$. Similarly, the Cauchy
problem (1),(B) has a continuous solution in a domain $\widetilde
D^+$.

In order to find the unknown coefficients in series (3), we
substitute expansions (3) into (1). After that, as usually, we
obtain a system of recurrence relations. We have the following
Cauchy problems for leading terms and terms in order $\varepsilon$
$$ \left\{\begin{array}{ll} u_{0t}+\lambda(u_0)u_{0x} =0, \
u_0|_{t=0}=\stackrel{o}{u}(x), x > 0 (< 0),
\\
\\ v_{0t}+\mu(u_0,v_0)v_{0x}=0, \ v_0|_{t=0}=\stackrel{o}{v}(x), x > 0 (< 0),
\end{array}\right. \eqno{(4.1)}$$
$$ \left\{\begin{array}{ll}
u_{1t}+\lambda(u_0)u_{1x}+\lambda_u(u_0)u_1u_{0x} =f(u_0,v_0), \
u_1|_{t=0}=0, x > 0 (< 0),
\\
\\ v_{1t}+\mu(u_0,v_0)v_{1x}+\mu_u(u_0,v_0)u_1v_{0x}+\mu_v(u_0,v_0)v_1v_{0x}
=g(u_0,v_0), \ v_0|_{t=0}=0, x > 0 (< 0).
\end{array}\right. \eqno{(4.2)}$$
Two various Cauchy problems for each system here are written down.
In the first case the initial data are taken on the semiaxis
$x>0$. In the second case the initial data are taken on the
semiaxis $x<0$. Thus, the asymptotic solution for the Cauchy
problem (1),(2) in the domains $\widetilde D^\pm$ have been
constructed. Moreover, in the same way, we can find all functions
$u_i, v_i$ in the domains $\widetilde D^\pm$.

The following properties take place (see [1]):
$D^-\subset\widetilde D^-$ and $D^+\subset\widetilde D^+$. Thus
the asymptotic solution for the Cauchy problem (1),(2) is found in
the domains $D^\pm$.
\section*{An asymptotics of shock lines
and an asymptotics of the solution in  domain $\widetilde D$.}
Let's remind, that the Hugoniot conditions follow from the
conservation laws of system (1). In this case, the Hugoniot
conditions looks as follows
$$ D^{\pm}\{u(t,s^{\pm}(t,\varepsilon),\varepsilon)- \widetilde
u(t,s^{\pm}(t,\varepsilon),\varepsilon)\}=
\Lambda(u(t,s^{\pm}(t,\varepsilon),\varepsilon))-
\Lambda(\widetilde u(t,s^{\pm}(t,\varepsilon),\varepsilon)),
\eqno{(5)}$$ $$ D^{\pm} \{ \Phi
(u(t,s^{\pm}(t,\varepsilon),\varepsilon),
v(t,s^{\pm}(t,\varepsilon),\varepsilon))- \Phi (\widetilde
u(t,s^{\pm}(t,\varepsilon),\varepsilon), \widetilde
v(t,s^{\pm}(t,\varepsilon),\varepsilon))\} = $$ $$ =\Psi
(u(t,s^{\pm}(t,\varepsilon),\varepsilon),
v(t,s^{\pm}(t,\varepsilon),\varepsilon))- \Psi (\widetilde
u(t,s^{\pm}(t,\varepsilon),\varepsilon), \widetilde
v(t,s^{\pm}(t,\varepsilon),\varepsilon)). \eqno{(6)}$$ Here $
\Lambda_u(u)=\lambda(u), \ D^\pm(t,\varepsilon) =
s_t^\pm(t,\varepsilon) $\footnote{The functions $\Phi$ and $\Psi$
are determined by a choice of the conservation laws of the system
(1).}. From the first conservation law we have the equations (5).
From the second we have the equations (6). The tilde notes that
the value of the discontinuous function is taken from the
$\widetilde D$ domain. As usually, we substitute the expansions
(2),(3) into the (6),(7) and use the asymptotic expansion for all
functions as $\varepsilon\rightarrow 0$. After that, collecting
terms in order $\varepsilon^0$ and $\varepsilon^1$ we obtain
$$D_0^\pm(u_0^\pm-\widetilde
u_0^\pm)=\Lambda(u_0^\pm)-\Lambda(\widetilde u_0^\pm), $$ $$
D_0^\pm\{\Phi(u_0^\pm,v_0^\pm)-\Phi(\widetilde u_0^\pm,\widetilde
v_0^\pm) \} = \Psi(u_0^\pm,v_0^\pm)-\Psi(\widetilde
u_0^\pm,\widetilde v_0^\pm),\eqno{(7.1\pm)}$$
 $$D_1^\pm(u_0^\pm-\widetilde
u_0^\pm)+D_0^\pm(u_1^\pm+u_{0x}^\pm s_1^\pm-\widetilde
u_1^\pm-\widetilde u_{0x}^\pm s_1^\pm ) =
\Lambda_u(u_0^\pm)(u_1^\pm+u_{0x}^\pm
s_1^\pm)-\Lambda_u(\widetilde u_0^\pm)(\widetilde
u_1^\pm+\widetilde u_{0x}^\pm s_1^\pm),$$
$$D_1^\pm\{\Phi(u_0^\pm,v_0^\pm)-\Phi(\widetilde
u_0^\pm,\widetilde v_0^\pm) \}+
D_0^\pm\{\Phi_u(u_0^\pm,v_0^\pm)(u_1^\pm+u_{0x}^\pm
s_1^\pm)+\Phi_v(u_0^\pm,v_0^\pm)(v_1^\pm+v_{0x}^\pm s_1^\pm)-$$
$$- \Phi_u(\widetilde u_0^\pm,\widetilde v_0^\pm)(\widetilde
u_1^\pm+\widetilde u_{0x}^\pm s_1^\pm)-\Phi_v(\widetilde
u_0^\pm,\widetilde v_0^\pm)(\widetilde v_1^\pm+\widetilde
v_{0x}^\pm s_1^\pm) \} = $$ $$ =
\Psi_u(u_0^\pm,v_0^\pm)(u_1^\pm+u_{0x}^\pm
s_1^\pm)+\Psi_v(u_0^\pm,v_0^\pm)(v_1^\pm+v_{0x}^\pm s_1^\pm)-$$
$$- \Psi_u(\widetilde u_0^\pm,\widetilde v_0^\pm)(\widetilde
u_1^\pm+\widetilde u_{0x}^\pm s_1^\pm)-\Psi_v(\widetilde
u_0^\pm,\widetilde v_0^\pm)(\widetilde v_1^\pm+\widetilde
v_{0x}^\pm s_1^\pm),\eqno{(7.2\pm)}$$ All function with tilde and
functions $s^\pm_i$, $D^\pm_i = s^\pm_{it}$ are unknown here.

Thus, from the Hugoniot conditions, we have four equations for six
unknown functions in every order of $\varepsilon $. Two additional
equations we can obtain from differential equations (1). These
equations are the same as the equations (4), when we replace both
$u$ and $v$ on $\widetilde v$ and $\widetilde v$ respectively.
Moreover, it is necessary to add the initial data to the
differential equations:
$$ \widetilde u_1|_{x=s^+_0(t)}=\widetilde u^+_1(t) \ \widetilde
v_1|_{x=s^-_0(t)}=\widetilde v^-_1(t).\eqno{(8.1)}$$. The leading
terms $\widetilde u_0, \widetilde v_0, s_0$ are the unperturbed
solution in domain $\widetilde D$.

Without loss of generality we can consider, that $v_0$ has a jump
on the curve $x=s^+_0(t)$, and $u_0$ is continuous function. Then
on the curve $x=s^-_0(t)$ the function $u_0$ has a jump and $v_0$
is continuous. Using this fact, we can determine all six unknown
functions. For example the procedure of a construction of $ u _ 1,
v _ 1, s ^ \pm _ 1 $ in domain $\widetilde D$ looks as follows

1) The function $u_0$ is continuous on the curve $x=s^+_0(t)$.
Hence we can solve the first equation in (7.2+) to determine
$\widetilde u_1^+$, because  $D_0^+, u_1^+$ are known.

2) To find the function $u_1$ inside of the domain $\widetilde D$
we must solve the Cauchy problem (4.2) for the first equation with
the initial data (8.1). In the same time, we determine the
function $\widetilde u_1^-$.

3) After that we can obtain the solution of the first equation in
(7.2-) to find $D_1^-$.

4)From the second equation in (7.2-) we can determine $\widetilde
v_1^-$ in the same way as an item 1.

5) To find the function $\widetilde v_1$ in the area $\widetilde
D$ and $\widetilde v_1^+$ respectively we must solve the Cauchy
problem (4.2) with the initial data (8.1).

6) To obtain $D_1+$ we solve the second equation in (7.2+).

This algorithm is convenient for determination of all coefficients
of the asymptotic expansion.

This work was supported by RFFI 00-01-00663, 00-15-96038 and INTAS
1068.

\end{document}